\documentclass[times,doublespace]{amsart}

\usepackage[mathscr]{eucal}
\usepackage{amssymb, amsmath,array, amscd}
\usepackage{enumerate}
\usepackage{hyperref}
\usepackage[all]{xy}

\def\rev{\usepackage[active]{srcltx}}
\def\p#1{\ensuremath{\mathbb{P}^{#1}}}
\def\eg{{\em e.g., }}
\def\ie{{\em i.e., }}
\def\C{\ensuremath{\mathbb C}}
\def\ve{\ensuremath{^\vee}}
\def\lar{\ensuremath{\,\longrightarrow\,}}
\def\be{\begin{equation}}\def\ee{\end{equation}}\def\ba#1{\begin{array}{#1}}\def\ea{\end{array}}\def\cl#1{\ensuremath{\mathcal{#1}}}
\def\st{\ensuremath{\,|\,}}
\def\pf{\noindent{\bf Proof. }}
\def\vez {\ensuremath{\times}}
\def\ra{\ensuremath{\,\rightarrow\,}}
\def\us{\ensuremath{^\star}}
\def\ls{\ensuremath{_\star}}
\def\ds#1{\displaystyle{#1}}
\def\M{{\cl  M}}
\def\D{{\cl  D}}
\def\na#1{\noalign{\vskip#1pt}}
\def\wt#1{\ensuremath{\widetilde{#1}}}

\newtheorem{coisa}{}[section]
\def\coi#1{\begin{coisa}{\bf
      #1.} \em}
\def\eco{\end{coisa}}
\def\cois#1{\begin{coisa}{\bf #1.}
\vspace{-10pt}
\eco
\addcontentsline{toc}{subsection}
{\numberline{}{#1}}
}

\def\runningheads#1#2{\markboth{\MakeTextUppercase{#1}}{\MakeTextUppercase{#2}}}

\def\rk{\ensuremath{\operatorname{rank}}}
\def\Hom{\ensuremath{\operatorname{Hom}}}
\def\cod{\ensuremath{\operatorname{cod}}}

\def\V{\ensuremath{\mathbb{V}}}
\def\ide#1{\ensuremath{\langle{#1}\rangle}}
\def\s#1{\ensuremath{\operatorname{Sym}_{#1}}}
\def\w{\ensuremath{\omega}}

\def\E{\ensuremath{\mathcal{E}}}
\def\schub{{\sc schubert}\,\cite{schub}}

\newcommand{\prin}{{\mathcal P}}

\newcommand{\de}{\partial}
\newcommand{\OX}{\mathscr{O}}

\newcommand{\T}{{\mathcal T}}
\newcommand{\Ker}{{\ensuremath{\mathrm{ker}}}}
\renewcommand{\ker}{{\ensuremath{\mathrm{ker}}}}

\newcommand{\Pe}{{\mathcal P}}

\rev

\begin{document}

\runningheads{V. Ferrer and I. Vainsencher}
{Singularities of Foliations}
\title{Degenerate singularities of one
  dimensional foliations}
\author{V. Ferrer
and I. Vainsencher}

\address{
IMPA --- Estrada Dona Castorina 110. 224620-320 Rio de
Janeiro Brazil. 
\newline\indent
ICEX--Depto. Matem\'{a}tica, UFMG -- Av. Antonio
Carlos 6627. 31270-901 Belo Horizonte Brazil. 
}

\thanks{The authors were partially supported by  CNPQ}


\begin{abstract}
We give formulas for the degrees of the spaces
of foliations in \p2 with a dicritical
singularity of prescribed order. Blowing up
such singularity induces, generically, a
foliation with all but finitely many leaves
transversal to the exceptional line; we 
also find the degree of the locus defined by
imposing a leaf of total contact with the
exceptional line.
\end{abstract}


\maketitle

\section*{Introduction}
\label{section Introduction}

Holomorphic foliations are an offspring of the
geometric theory of polynomial differential 
equations. Following the trend of
many branches in Mathematics,
interest has migrated to global aspects.
Instead of focusing on just
one curve, or surface, or
metric, or differential equation, try and study
their family in a suitable parameter
space. The geometry within the parameter space
of the family acquires relevance.  For
instance, the family of hypersurfaces of a
given degree correspond to points in a 
projective space;  geometric conditions on
hypersurfaces usually correspond to interesting
subvarieties in the parameter space, \eg the
discriminant.  Hilbert schemes have their
counterpart in the theory of polynomial
differential equations, to wit, the spaces of
foliations.

Holomorphic foliations of degree $d$ on the
complex projective plane \p2 are defined by
nonzero twisted 1--forms, $\omega=\sum
a_idz_i$, with homogeneous polynomials
$a_i(z_0,z_1,z_2)$ of degree ${d+1}$, up to
scalar multiples, satisfying $\sum a_iz_i=0$.
The parameter space of foliations of degree $d$
is a projective space \p N (cf.\,\ref{H0}).

The purpose of this work is to compute the
dimensions and degrees of the subvarieties of
\p{N} corresponding to foliations displaying
certain degenerate singularities. Given an
integer $k\geq2$ we study the locus,
$M_k\subset\p N$, of foliations with a
singularity of order $\geq k$.  These are
foliations defined in local coordinates by a
holomorphic 1--form that can be written as
$\omega=a_kdx+b_kdy+\rm higher\ order\ terms$,
with $a_k(x,y),b_k(x,y)$ homogeneous
polynomials of degree $k$. It turns out that
$M_k$ is the birational image of an explicit
projective bundle over \p2. This enables us to
find a formula for the degree of $M_k$.

Another interesting type of non-generic
foliation presents a so called {\em dicritical}
singularity of order $k$: require $a_kx+b_ky$
to vanish.  This defines a closed subset
$D_k\subset M_k$.  A characteristic feature of
a foliation with a dicritical singularity is
the fact that all but finitely many leaves of
the foliation induced on the blowup at the
singular point are transversal to the
exceptional line.  We say a foliation with such
a singularity has the property {\em of maximal
contact} ({\em MC}\, for short) if some leaf of
the induced foliation has a contact of order
$k$ with the exceptional line. Thus we may
consider the subvariety $C_k\subset D_k$
consisting of dicritical foliations with {\em
MC}.

Requiring a leaf of a foliation to be tangent
to a line at a given point defines a hyperplane
in the parameter space \p N. Therefore,  the
degree of each of the loci 
$C_k\subset{}D_k\subset{}M_k\subset\p N$ can be rephrased
loosely as  the number of foliations
with a singularity of the chosen type and
further tangent to the appropriate number of
flags (point,\,line) in \p2.  It turns out that
the degrees of $C_k,D_k,M_k$ are expressed as
explicit polynomials in $k,d$.

This fits into the tradition of classical
enumerative geometry: answers to questions such
as determining the number of plane algebraic
curves that have singularities of prescribed
orders, besides passing through an appropriate
number of points in general position, are often
given by ``node'' polynomials.
There is also a wealth of results and
conjectures on generating functions for
counting suitably singular members of linear
systems of curves on surfaces,
cf.\,G\"otsche\,\cite{goetsche}, Kleiman and
Piene\,\cite{KP}. We hope similar results can
be formulated in the setting of foliations.

\section{The space of foliations}

The main reference for this matterial is
Jouanolou\,\cite{Joua}. A projective 1-form of
degree $d$ in $\p 2$ is a global section of
$\Omega^1_{\p 2}\otimes \OX_{\p 2}(d+2)$, for
some $d\geq 0$.  

We denote by $S_d$ the space 
$ 
{\rm Sym}_d(\C^{3})\ve
$ of
homogeneous polynomials of degree $d$ in the
variables $z_0,z_1,z_2$.  We write
$\partial_i=\partial/ \partial z_i $, thought
of as a vector field basis for $\C^3$. The dual
basis will also be written as $dz_0,dz_1,dz_2$
whenever we think of differential forms.
Recalling Euler sequence
\begin{equation}\label{Eulerdual}
0\to \Omega^1_{\p 2}(d+2)\to \OX_{\p
  2}(d+1)\otimes S_{1}\to \OX_{\p 2}(d+2)\to 0
\end{equation} 
and taking global sections
we get
 the exact sequence
\begin{equation*}
0\to H^0(\p2,\Omega^1_{\p2}(d+2))\lar S_{d+1}\otimes S_{1}\stackrel{\iota_{R}}\lar S_{d+2}\to 0 
\end{equation*} 
where $\iota_{R}(\sum a_{i}dz_{i})=\sum a_{i}z_{i}$ is the
contraction by the radial vector field.
Thus a 1-form $\w\in
H^0(\p2,\Omega^1_{\p2}(d+2))$ can be 
written  in homogeneous coordinates as 
$$\w=a_{0}dz_{0}+a_{1}dz_{1} +a_{2}dz_{2}$$ where the $a_{i}$'s are
homogeneous polynomials of degree $d+1$ satisfying  
$$
a_{0}z_{0}+a_{1}z_{1} +a_{2}z_{2}=0.
$$ 
The space of foliations of degree $d$ in \p2 is
the projective space
\be\label{H0}
\ba c
\p N
=\p{}(H^0(\Omega^1_{\p2}(d+2)))
\ea\ee
of dimension
\\\centerline{$
N=3\binom{d+2}2-\binom{d+1}2-1
=d^2+4d+2
$.}  
We have $\Omega^1_{\p
2}\otimes \OX_{\p 2}(d+2) =\Hom(\T\p2,\OX_{\p
2}(d+2)) $. Any nonzero
\\\centerline{$
\omega:\T\p2\lar\OX_{\p2} (d+2)
$} 
induces a (singular, integrable) distribution
of dimension one subspaces given by $p\mapsto
\Ker \w_{p}$.  A nonzero multiple of \w\ yields
the same distribution.

\coi{Singularities}\label{sec:sing}
The singular scheme of \w\ 
is defined by the ideal sheaf
image of 
$
\omega\otimes \OX_{\p2}(-d-2):
\T\p2\otimes\OX_{\p2}(-d-2)\ra
\OX_{\p 2}$. 
If finite, it consists of
$\int_{\p2}c_2(\Omega^1_{\p 2}(d+2))= d^2+d+1$ 
\ points. 
In local coordinates, 
say around $0=[0,0,1]\in\p2$, 
 writing $\w=adx+bdy$ the
singular scheme of $\w$ is given by the ideal
$\ide{a,b}$.  We say $0$ is a
nondegenerate singularity if the jacobian
determinant $|\partial(a,b)/\de(x,y)|$
is nonzero.
A generic 1-form has only nondegenerate 
(hence isolated) singularities;
see~\cite[p.\,87]{Joua}.
The order of the singularity $0$ is 
$\nu_{0}(\w)=\min\{{\rm order}_{0}(a)
,{\rm order}_{0}(b)\}.
$
It can easily be checked that this is
independent of the choice of coordinates. 
\eco
In
fact, if $\cl I$ is the ideal sheaf of the
singular scheme, then for each $p\in\p2$ there
is a unique nonnegative integer  $k$ such that
the stalk at $p$ 
satisfies $\cl I_p\subset\mathfrak m_p^k$
and $\cl I_p\not\subset\mathfrak m_p^{k+1}$,
where $\mathfrak m$ denotes the ideal sheaf of
$p$. Thus, setting $\cl
E=\Omega^1_{\p2}\otimes\OX_{\p2}(d+2)  
$, we see that
the order of the singularity $p$ is
at least $k$ iff the image of $\w$ in the
quotient $\cl E/\mathfrak m_p^k\cl E$  is zero.

\coi{Jet bundles}\label{sec:jet}
The preceeding discussion entices us to
 recall the notion of jet bundles associated to
 a vector bundle, cf.\,\cite[16.7]{EGA},
\cite{P}. 
Let $\E$ be a vector bundle over a smooth
projective variety $X$. For $k\geq 0$ the
$k$th-jet bundle associated to $\E$, denoted
$\Pe^k(\E)$, is a fiber bundle over $X$ with
fiber over  $x\in X$ given by
\\\centerline{$
\Pe^k(\E)_{x}=(\OX_{X}/m_{x}^{k+1})\otimes
\E_{x}
$} 
where $m_{x}$ is the maximal ideal of the point $x$.
\eco
For each $k\geq 0$ we have exact sequences
\begin{equation}\label{jetn}
0\to \s{k+1}\Omega^1_{X}\otimes \E\to \Pe^{k+1}(\E)\to \Pe^k(\E)\to 0.
\end{equation}

Consider the evaluation map 
$$ev:X\times H^0(X,\E)\to \E
$$ given by
$ev(x,s)=(x,s(x))$. 
The map $ev$  lifts to natural maps fitting
into the following commutative diagram:
\be\label{evk}
\xymatrix
{\,X\times H^0(X,\E)\ar@{->}[r]
^{\hskip.71cm
ev_{k}}
& \Pe^{k}(\E)
\ \ar@{->}[d] 
\\
&\Pe^{k-1}(\E)
\ \ar@{<-}[ul] 
^{ev_{k-1}} .
}
\ee
We think of $ev_k(s)$ as the Taylor
expansion of $s$ truncated at order $k+1$.  
We include for the reader's convenience the
following 
\coi{Lemma of global generation}
\label{sec:global}\em
Notation as in {\rm(\ref{evk})} above,
given $k$,
replacing \cl E\ by a
sufficiently high twist $\cl E\otimes\cl L^m$
by an ample line bundle \cl L, we have that
\\{\rm(i)} the map $ev_k$
is surjective.
\\{\rm(ii)}
Set
 $W_x=\{s\in{}H^0(X,\E)\st(ev_{k-1})_x(s)=0\}$.
Then
\\\centerline{$
(ev_k)_x(W_x)=(\s{k}\Omega^1_{X}\otimes \E)_x
$.}

\eco

\pf
Let \cl J\ be the ideal of the diagonal of
$X\vez{}X$. Consider the projection maps $p_i:X
\vez{}X\ra{}X,\,i=1,2$. We have the exact
sequence of sheaves over $X\vez{}X$,
\\\centerline{$
\xymatrix@R1pc
{
\cl J^{k+1} \   \ar@{^(->}[d] \ 
\ar@{=}[r]\ &
\cl J^{k+1}
\ar@{^(->}[d]
&
\\
\cl J^{k} \   \ar@{->>}[d]
\ar@{^(->}[r]&
\OX\ar@{->>}[d]
\ar@{->>}[r]&
\OX/\cl J^{k}
\ar@{=}[d]
\\
\cl J^{k}/\cl J^{k+1}\ 
\ar@{^(->}[r]&
\OX/\cl J^{k+1}
\ar@{->>}[r]&
\OX/\cl J^{k}.
}
$
}
Tensoring by $p_2\us\cl E\otimes\cl L^m$ and
taking $(p_1)\ls$ yields
\\\centerline{$
\xymatrix@R1.1pc@C.7pc
{&
(p_1)\ls((\cl O/\cl J^{k+1})\otimes\cl E\otimes\cl L^m)
\ar@{=}[d]&\\
\s{k}\Omega^1_X
\otimes\cl E\otimes\cl L^m\ 
\ar@{>->}[r]&
\prin^k(\cl E\otimes\cl L^m)
\ar@{->>}[r]&
\prin^{k-1}(\cl E\otimes\cl L^m)
.
}
$
}
Surjectivity of $ev_k$ follows upon killing
$(R^1p_1)\ls(\cl J^{k+1}\cl E\otimes\cl L^m)$.
Similarly, (ii) follows from the identification
$W_x=(p_1)\ls(\cl J^{k}\cl E\otimes\cl
L^m)_x$. \qed

\coi{Singularities of order $k$}
We apply the previous lemma to $\cl
E=\Omega^1_{\p 2},\cl L=\OX_{\p 2}(1)$. 
In order to simplify the notation we set for
short in the sequel
\eco
\centerline{$\Omega=\Omega^1_{\p2}\ \text{ and }
\ V:=H^0(\Omega(d+2))
$.}
\coi{Remark}\label{surj}
Fix $k\leq d+1$. It follows from the explicit
calculation of $H^i(\p2,\Omega(d+2))$
that the conclusions of the previous lemma hold
for $m=d+2$.
Hence 
\\\centerline{$\ba l
\text{
{\bf(i)} $ev_k$ is surjective for
$k\leq d+1$ and} 
\\
\text{
{\bf(ii)} 
$ev_{k}(\ker(ev_{k-1}))=
\s{k}\Omega\otimes\Omega (d+2).
$}\ea$}
\eco
 
We describe the locus
$M_{k}\subset \p {N}$ of foliations of given
degree $d$ that have some singularity of order
$\geq k$.

\begin{coisa}{\bf Proposition.}\label{degMk}
For $1\leq k\leq d+1$, denote by
$$M_{k}=\{[\w]\in\p N\mid [\w]\,\text{ \rm 
has a singularity of order at least }\,\,
k\}.
$$
Then we have
\\\centerline{$
\cod_{\p N}
M_{k}=k(k+1)-2
$}
\\ and 
\\\centerline{$
\deg(M_{k})=\ds\int_{\!\p2}\,
c_{2}(\prin^{k-1}(\Omega(d+2)))
.$}

\end{coisa}
\begin{proof}

Define  
\\\centerline{
$\M_{k}=\Ker\left(
ev_{k-1}:\p2\vez{}V\ra\prin^{k-1}(\Omega (d+2))\right)
$.}
In view of the previous remark, we see that
 $\M_k$ is a vector subbundle of $V$ of
co-rank equal to $\rk\prin^{k-1}(\Omega 
(d+2))$.
By construction, the projective bundle
associated to $\M_{k}$ is the incidence variety,
$$\p{}(\M_{k})=\{(p,[\w])\in
 \p 2\times \p N\mid p\,\text{ is a
  singularity of }\,[\w]\,
 \text{ and }\, \nu_{p}(\w) \geq k\}.
$$
Let $q:\p{}(\M_{k})\to \p N$
denote the projection in the second factor. 
We have $M_k=q(\p{}(\M_{k}))$.
It is easy to check  that 
$q$ is generically injective. It follows that
 $\deg(M_{k})=\int s_{2}(\M_{k})\cap [\p 2]$.
Since
$s_{2}(\M_{k})=c_2(\prin^{k-1}(\Omega 
(d+2)))$, 
the assertions now follow from (\ref{jetn}).
\end{proof}

%
%
%
%
%
Using  the proposition  
we may now derive an explicit formula
for  the degree of $M_{k}\subset\p N$. 
See also the
script in \S\,\ref{script}.
We find

\begin{coisa}{\bf Corollary.}\label{cor2.1}
The degree of $M_{k}$ is given by 
$$\frac{1}{2}k(k+1)\left[
(k^2+k-1)\left(d^2-(2k-3)
d\right)+
\frac{1}{4}(4k^4-8k^3-7k^2+21k-6)\right].
$$
~\hfill\qed
\end{coisa}

 
\section{Dicritical singularities}

If $\w\in H^0(\p 2,\Omega 
(d+2))$ and $p$ is a
singularity of $\w$, we say that $p$ is
dicritical if the local expression of $\w$ is 
\\\centerline{$\w_{p}=a_{k}dx+b_{k}dy+h.o.t
$} with
$a_{k}x+b_{k}y=0$.
In the case $k=1$, we say that $p$ is a {\em
  radial singularity}.  

Observe that this  condition is equivalent to 
$$\w_{p}=f(x,y)(ydx-xdy)+h.o.t$$ for some
homogeneous  polynomial $f$ of degree $k-1$.

The main result of this section is the
following.

\begin{coisa}\label{Dk}
{\bf Proposition.}
For all $1\leq k\leq d$ there exists a
subbundle $\D_{k}$ of $\M_k\ra\p 2$ such that 
\\$
{\rm(i)}\ \p{}(\D_{k})=\{(p,[\w])
\in  \p 2\times \p N
\mid p\,\text{ is a
  dicritical singularity
  of }\,[\w]\,\text{with}\,\nu_{p}(\w)\geq
k\}.
$
\\{\rm(ii)} 
Set   $D_{k}=q(\p{}(\D_{k}))$.  
Then the
codimension of  $D_{k}$ in \p N is $k(k+2)$ and
\\{\rm(iii)} 
the
 degree of  $D_{k}$ is the coefficient  of  the
 degree two part of 
\\\centerline{$
c(\prin^{k-1}(\Omega (d+2)))c(\s{k+1}
\Omega \otimes\OX_{\p 2}(d+2)).
$}
 
\end{coisa}

\coi{\bf Remark}\label{multil}
Before proceeding to the proof of the
proposition, we explain an
 invariant way of expressing the
condition that a singularity be dicritical.
Suppose that  $\E$ is a vector bundle of rank $2$.
Then for all $k\geq 1$ we have the following  exact sequence 
(\eg see \cite[Appendix 2 A2.6.1.]{Eis}),
$$0\to\stackrel2\wedge\E\otimes \s{k-1}\E\to
\s{k}\E\otimes \E\stackrel{P_{k}}\to \s{k+1}\E\to
0,
$$
where the first map is given by
$$(a\wedge b\otimes c)\mapsto (ac\otimes
b)-(bc\otimes a)
$$
and the second by multiplication, \ie
$$a\otimes b\mapsto ab.
$$
Say $x,y$ form a local basis for \cl E. Then
for  $a_k,b_k\in\s k\cl E$, we have
that $a_kx+b_ky=0$ in \s{k+1}\cl E iff there is
some $c\in\s{k-1}\cl E$ such that
$a_k\otimes{}x+b_k\otimes{}y
$ is equal to 
the image of $x\wedge{}y\otimes c$, to wit,
$xc\otimes{}y
-yc\otimes{}x
$.
\eco
\coi{Construction of $\cl D_k$}
We have the following diagram,
$$
\xymatrix
{
&& \s{k}\Omega\otimes\Omega (d+2)\ 
\vphantom{\ba c1\\.\ea}
\ar@{>->}[d]
\\
&&\hskip.31cm
\prin^{k}(\Omega (d+2))\ \ar@{->>}[d]
\ \ar@{<<-}@/_0.6cm/[ld]^{
\hskip-1cm
ev_k\ }& 
\\
\cl M_k \ \ar@{>->}[r]
\ \ar@{->>}@/^0.6cm/[uurr]^{\hskip-1.851cm
J_k}
&\ V\ \ar@{->>}[r]& \prin^{k-1}(\Omega (d+2))
}
$$
The map $J_{k}$ defined in the previous diagram
is surjective in view of remark\,\ref{surj}.
We obtain the surjective map
\\\centerline{$
\xymatrix
{
\M_{k}
\ \ar@{->>}@/_0.6cm/[rr]_{T_{k}}
\hskip.1cm\ \ar@{->}[r]^{\hskip-.9841cmJ_{k}}
&
\s{k}\Omega \otimes\Omega (d+2)
\ \ar@{->}[r]^{P_{k}}
&\hskip.21cm
\s{k+1}\Omega(d+2).
}
$}
\eco
Explicitly, on the fiber over
 $p\in\p 2$ the map is as follows:  
$$T_{k}(p,\w)=(p,a_{k}x+b_{k}y)$$ where 
\\\centerline{$
\w_{p}=a_{k}dx+b_{k}dy+h.o.t.
$} is the local expression of $\w$ in a neighborhood of $p$.
Set 
\be
\label{defDk}
\D_{k}:=\ker\left(
\vphantom{I^I_I}\right.
\xymatrix
{
\M_k\ \ar@{->>}[r]^{
\hskip-1.043cm
T_{k}}& 
\s{k+1}\Omega (d+2)
}
\left.\vphantom{I^I_I}
\right)
\ee

 Thus $\D_{k}$ is a vector bundle of rank 
$=\rk(\M_{k})-(k+2)$.
Recalling  \ref{multil}, we see that
 the projective bundle associated to 
 $\D_{k}$ is the incidence variety,
\\\centerline{$
\p{}(\D_{k})=\{(p,[\w])
\in  \p 2\times \p N
\mid p\,\text{ is a
  dicritical singularity
  of }\,[\w]\,\text{ with }\,\nu_{p}(\w)\geq
k\}.
$}
It can be shown
 that $q$ is generically injective. Hence
the degree of $D_{k}$ is given by 
$\int s_{2}(\D_{k})\cap [\p 2].
$
This finishes the proof of the
 proposition\,\ref{Dk}.
\qed

A formula for the degree of $D_k$ can be made
explicit.

\begin{coisa}{\bf Corollary.}\label{degDk}
The degree of $D_{k}$ is given by
$$(k+1)^2
\left
[\frac{1}{2}(k^4+k^2-2k+2)-(k^3+k^2+k-1)d+
\frac{1}{2}(k^2+2k+2)(k+1)^2d^2
\right].
$$
~\hfill\qed
\end{coisa}

\begin{coisa}{\bf Remarks.}\label{obs1}
\em {\bf(i)}
 We
have by construction  the following diagram:  
\[
\xymatrix
{
          &  &\s{k-1}\Omega\otimes \wedge^2\Omega\otimes \OX_{\p 2}(d+2)\ar@{->}[d]\\
\D_{k} \ar@{->}[r]
\ \ar@{->}@/^0.6cm/[urr]_{
\hskip-1.8cm
d_{k}}
\ \ar@{->}[drr]_{0\ }&\M_{k}\ar@{->}[r]^
{
\hskip-1cm
J_{k}}\ar@{->}[dr]^{T_{k}}& \s{k}\Omega\otimes \Omega(d+2) \ar@{->}[d]^{P_{k}}\\
          &  &  \ \ 
\s{k+1}\Omega\otimes \OX_{\p 2}(d+2)\ \ 
}
\]

By definition  of $\D_{k}$  we obtain a map $$d_{k}:\D_{k}\to
\s{k-1}\Omega\otimes \wedge^2\Omega\otimes
\OX_{\p 2}(d+2)
$$ 
given in the fibers by $d_{k}(p,\w)= f(x,y)dx\wedge dy$ where $f$ is a
polynomial of degree $k-1$.  

\medskip
{\bf(ii)}
 In the case $k=1$ we have
  $$\w=\lambda(ydx-xdy)+h.o.t.
$$
with $\lambda\in \C$ \ie {\em a radial} 
singularity. 
Thus \ref{Dk} and 
 \ref{degDk} give 
formulas for the codimension and degree of the
space of  foliations with a radial
singularity:
$$\left\{
\ba c
\cod_{\p N} D_{1}=3 \\\na{7}
\deg D_{1}=10d^2-8d+4.
\ea
\right .
$$

{\bf(iii)}
 In the case $k=d+1$  the map
 $J_{d+1}:\M_{d+1}\to \s{d+1}\Omega^1_{\p
   2}\otimes \Omega 
(d+2)$ 
is no longer surjective: its image  is 
$\s{d}\Omega \otimes \wedge^2\Omega^1_{\p
  2}
\otimes \OX_{\p 2}(d+2)$. Indeed, suppose that $\w$ 
is a form of degree $d+1$ which
has $p$ as  singularity of order
$d+1$. Then a local expression of $\w$ is 
\\\centerline{$
\w_{p}=a_{d+1}dx+b_{d+1}dy
$}
 but this form defines a projective
form of degree $d+1$ in 
$\p 2$ if and only if $$a_{d+1}x+b_{d+1}y=0$$ i.e. if $p$ is a  dicritical  singularity.
Therefore we can
write $$\w_{p}=f(x,y)(ydx-xdy)$$ for some
homogeneous polynomial  $f$  of degree $d$, i.e. $\w_{p}\in \s{d}\Omega^1_{p}\otimes \wedge^2\Omega^1_{p}$.
Hence 
\\\centerline{$
\xymatrix
{T_{d+1}:\M_{d+1}\ \ar@{->}[r]^{
\hskip-.61cm
J_{d+1}}
& \s{d+1}\Omega \otimes\Omega (d+2)
\ \ar@{->}[r]^{\!\!P_{d+1}}&
\s{d+2}\Omega \otimes\OX_{\p 2}(d+2)}
$}  
is the zero map.  
This shows that $\M_{d+1}=\D_{d+1}$ \ie, for a
foliation of degree $d$ a singularity of order
$d+1$ is automatically dicritical.    

\end{coisa}

\section{Maximal contact}

Consider a degree $d$ form  $\w\in  H^{0}(\p 2,\Omega (d+2))$ with a dicritical singularity at $p=[0:0:1]$ of order  $k\geq 2$. Denote by  
$$\pi:\widetilde{\C}^2\to \C^2$$ the blowup of $\C^2$ at $p$.   
 Write $\w$ in
local coordinates $(x,y)$ around $0\in\C^2$ as
 $$\w=\sum_{j=k}^{d+1} a_{j}dx+b_{j}dy$$
where $a_{j},b_{j}$ are homogeneous polynomials of degree $j$.
The blowup of $\C^2$ at $(0,0)$  
$$\widetilde{\C}^2=\{(x,y),[s:t]\mid tx=sy\}\subset \C^2\times \p 1$$ is
covered by the usual two charts 
$$V_{0}=\{((x,y),[1:t])\mid tx=y\}\simeq
\{(x,t)\mid t,x\in \C\},
$$ 
$$V_{1}=\{((x,y),[s:1])\mid x=sy\}\simeq
\{(s,y)\mid s,y\in \C\}.
$$
Over $V_{0}$ we have $dy=tdx+xdt$. Thus
\begin{align*}
\pi^*(\w)(x,t)=\sum_{j=k}^{d+1} a_{j}(x,tx)dx+b_{j}(x,tx)(tdx+xdt)=\\
\sum_{j=k}^{d+1} (a_{j}(x,tx)+tb_{j}(x,tx))dx+xb_{j}(x,tx)dt=\\
x^k\sum_{j=k}^{d+1}x^{j-k} 
\left[
(a_{j}(1,t)+tb_{j}(1,t))dx+xb_{j}(1,t)dt\right].\\
\end{align*}
Since $p=(0,0)$ is a dicritical singularity we have
$$\w=f(x,y)(ydx-xdy)+h.o.t$$ where $f$
is a polynomial of degree $k-1$; so
$a_{k}=yf(x,y)$, $b_{k}=-xf(x,y)$.

Hence we may write
\begin{align*}
\pi^*(\w)(x,t)=x^k(xb_{k}(1,t)dt+x\alpha)=x^{k+1}(-f(1,t)dt+\alpha)
\end{align*}
where $\alpha$ is a 1-form.
The strict transform of $\w$ is
\begin{equation}\label{blow}
\widetilde{\w}=-f(1,t)dt+(a_{k+1}+tb_{k+1})dx+x\alpha_{1}
\end{equation}  
for some 1-form $\alpha_{1}$.
%
%
Over $V_{0}$ the exceptional divisor is given by
$x=0$, and  by  (\ref{blow}) we have 
$$\widetilde{\w}\wedge dx=-f(1,t)dt\wedge
dx+x\alpha_{1}\wedge dx.
$$
The leaves of 
$\wt\w$ passing
through each point $(0,t_{0})$ with
$f(1,t_{0})\neq 0$ are transverse to the
exceptional divisor. On the other hand, the
points $(0,t_{0})$ such that $f(1,t_{0})=0$
but aren't singularities of $\widetilde{\w}$
are exactly the points of tangency of
$\widetilde{\w}$ with the exceptional divisor.

Next, we study the relationship between the
multiplicity of $t_{0}$ as a
zero of $f(1,t)$ and the order of tangency of
the leaf of $\widetilde{\w}$ 
with the exceptional divisor at $(0,t_{0})$.

\begin{coisa}{\bf Lemma.}
The intersection multiplicity of a leaf of 
$\wt\w$ with the exceptional divisor at a
point $(0,t_{0})$ is the multiplicity of $t_{0}$ as zero of $f(1,t)$ plus one.
\end{coisa}

\begin{proof}
We may assume  $t_{0}=0$.  By (\ref{blow})
we have that $\widetilde{\w}$ has the following
form
$$\widetilde{\w}=(-f(1,t)+xF(x,t))dt+(g(t)+xG(x,t))dx$$
with $f(1,0)=0$. Observe that $g(0)\neq 0$
because we are assuming that $p:=(0,0)$ is a
 nonsingular point of $\widetilde{\w}$.  Let
$h(x,y)=0$ be a local equation for a leaf of
$\widetilde{\w}$ through $(0,0)$, where $h$ is
a non constant holomorphic function.  We have
\\\centerline{$(\widetilde{\w}\wedge dh) (p)
=-g(0)\frac{\de h}{\de t}(p)dt\wedge dx=0
$.}  Hence $\frac{\de
h}{\de t}(p)=0\neq \frac{\de h}{\de x}(p)$.  Therefore, we
can find a local analytic
parameterization of $h=0$ of the form $y=t,x=\gamma(t)$
defined
in a neighborhood of $t=0$ such that
\\\centerline{$
\left\{
\ba l
\gamma(0)=0\\\na7
\gamma'(0)=
\frac{\de{}h}{\de{}t}(p)
\big/
\frac{\de h}{\de x}(p)
=0
\ea
\right.
$}

\medskip
Since  $(\gamma(t),t)$ parameterizes a leaf
of $\widetilde{\w}$ we find
that 
$$f(1,t)+\gamma(t)F(\gamma(t),t)
+\gamma'(t)(g(t)+\gamma(t)G(\gamma(t),t))\equiv
0.
$$
Hence, repeatedly 
differentiating with respect to $t$  yields
$$\gamma^{(j)}(0)=0\ \forall\, 
j\leq r \implies
\gamma^{(r+1)}(0)= \frac{\frac{\de^r f}{\de^r
t}(1,0)}{g(0)}.
$$
Now, $\gamma$ has intersection multiplicity $n$
with $x=0$ at $(0,0)$ if the
 first non-vanishing derivative of
$x(\gamma(t),t)=\gamma(t)$ at $0$ is precisely
$n$.  Thus the intersection multiplicity
 of $h=0$ with $x=0$ is $n$ if and
only if \ $t=0$ \ is a zero of order $n-1$ of
$f(1,t)$.
\end{proof}

\medskip

From the above Lemma we have that if $\w\in
D_{k}$, then the order of tangency of the
leaves of $\w$ with the exceptional divisor
is $\leq k$, and is equal to $k$ precisely in 
the case that $f=l^{k-1}$ where $l$ is a
polynomial of degree one.


\coi{Degree of the {\em MC} locus} Recall that we say
that a form $\w$ has the {\em MC} property if
it has a dicritical singularity $p$ of order
$k$ such that the strict transform of $\w$
under the blowup of $p$ has a leaf with maximal
order of contact with the exceptional divisor
of the blowup.

 \eco
Consider a form with a dicritical singularity
of order $k$, 
$$
\w=f(x,y)(ydx-xdy)+h.o.t.
$$ 
(i.e. $f$ is a polynomial of degree
$k-1$). Then $\w$ has the  \emph{MC}
property
if and only if $f(1,t)=(t-t_{0})^{k-1}$ or
$f(s,1)=(s-s_{0})^{k-1}$, i.e.  
$$f(x,y)=(ax+by)^{k-1}\,\text{ for some }\,
a,b\in \C.
$$
\indent{}
Therefore we can parameterize the set of forms
that has the \emph{MC} property as follows. The Veronese-type map $\Omega \to \s{k-1}\Omega $
induces an embedding 
$$\mathrm{v}_{k}:\p{}(\Omega \otimes\wedge^2\Omega \otimes
\OX_{\p 2}(d+2))\to
\p{}(\s{k-1}\Omega \otimes\wedge^2\Omega \otimes
\OX_{\p 2}(d+2))
$$
which locally is given by 
$\mathrm{v}_{k}(p,l\otimes (ydx-xdy))=(p,l^{k-1}\otimes(ydx-xdy))$.
In order to simplify the notation set 
$$\E:=\Omega \otimes\wedge^2\Omega \otimes
\OX_{\p 2}(d+2)$$ and
$$\E_{k}:=\s{k-1}\Omega \otimes\wedge^2\Omega \otimes
\OX_{\p 2}(d+2).
$$

Define
$$\V_{k}:=\mathrm{v}_{k}(\p{}(\cl E))\subset
\p{}(\E_{k})$$

\begin{coisa}{\bf Lemma.}\label{alpha} 
The codimension of $\V_{k}$ in
$\p{}(\E_{k})$ is $k-2$
and  its cycle class is
 $$[\V_{k}]=u H_{k}^{k-2}+v hH_{k}^{k-3}+w
h^2H_{k}^{k-4}\cap \p{}(\E_{k})$$
where  $H_{k}\,(resp. \ h)$ denotes
 the relative hyperplane class of
$\p{}(\E_{k})\,(resp.\ \p2)$  and
$$\ba cu=(k-1),
\\
v=-\frac{1}{2}(k-1)(k-2)(3k+2d-5),
\\
w=\frac{1}{8}(k-2)(k-1)^2(9k^2-47k+12kd-60d+72+12d^2).
\ea
$$
\end{coisa}

\begin{proof}

It is clear that 
$\cod(\V_{k})=
\dim\p{}(\E_{k})-\dim
\p{}(\E)=
k-2$.
Recalling that the Chow ring
$A_{*}(\p{}(\E_{k}))$ is
generated by $H_{k}$ and $h$ 
(see~\cite[Thm.\,3.3.,\,p.\,64]{Fulton}) 
we can express
\begin{equation}\label{V}
[\V_{k}]=u H_{k}^{k-2}+v hH_{k}^{k-3}+w
h^2H_{k}^{k-4}\cap \p{}(\E_{k}).
\end{equation}

With this notation, the relative hyperplane class of
$\p{}(\E)$ is $H_{2}$, and
we have $\mathrm{v}_{k}^*(H_{k})=(k-1)H_{2}$.
Consider the following diagram:

\[
\xymatrix 
{
\p{}(\E)\ar@{->}[dr]_{\rho}\ar@{->}[r]^{\mathrm{v}_{k}}  & \p{}(\E_{k})\ar@{->}[d]^{\pi}\\
 & \p 2  }
\]
To find the coefficient $u$  we multiply by
$h^2H_{k}$ both sides of (\ref{V}) to obtain:  
$$h^2H_{k}\cap \mathrm{v}_{k*}(\p{}(\E))=u
h^2H_{k}^{k-1}\cap\p{}(\E_{k}).$$
By  the projection formula we have:
$$h^2 \mathrm{v}_{k*}((k-1)H_{2}\cap \p{}(\E))=u
h^2H_{k}^{k-1}\cap\p{}(\E_{k}).
$$
Applying $\pi_{*}$ to this last equation we find
\begin{align*}
&(k-1)h^2\cap\rho_{*}(H_{2}\cap\rho^*\p 2)=u
  h^2\pi_{*}(H_{k}^{k-1}\cap\pi^*\p 2)\\
&(k-1)h^2\cap s_{0}(\E)\cap\p 2=u h^2\cap s_{0}(\E_{k})\cap
  \p 2.
\end{align*}
Hence $u=k-1$. Next, multiplying (\ref{V}) by 
$hH_{k}^2$  we get   
$$hH_{k}^2\cap \mathrm{v}_{k*}(\p{}(\E))=u
hH_{k}^{k}\cap\p{}(\E_{k})+v
h^2H_{k}^{k-1}\cap\p{}(\E_{k}).
$$
Hence the  projection formula gives us
$$h(k-1)^2\cap \mathrm{v}_{k*}(H_{2}^2\cap \p{}(\E))=u
hH_{k}^{k}\cap\p{}(\E_{k})+v
h^2H_{k}^{k-1}\cap\p{}(\E_{k}).
$$
Applying $\pi_{*}$ we obtain
\begin{align*}
&h(k-1)^2\rho_{*}(H_{2}^2\cap \p{}(\E))=u
h\pi_{*}(H_{k}^{k}\cap\pi^*\p 2)+v h^2\pi_{*}(H_{k}^{k-1}\cap\pi^*\p 2),\\
&(k-1)^2hs_{1}(\E)\cap \p 2=((k-1)hs_{1}(\E_{k})+v
h^2s_{0}(\E_{k}))\cap \p 2,
\end{align*}
hence $v=(k-1)^2s_{1}(\E)-(k-1)s_{1}(\E_{k})$.

Similarly we obtain 
$w=[(k-1)^3s_{2}(\E)-u
s_{2}(\E_{k})-v h s_{1}(\E_{k})]\cap
[\p 2]$.

The lemma follows from the calculation of the
Segre
classes of $\E$ and $\E_{k}$. Observe that
$\wedge^2\Omega =\OX_{\p 2}(-3)$, so that
 $\E=\Omega (d-1)$
and $\E_{k}=\s{k-1}\Omega \otimes \OX_{\p 2}(d-1)$. These classes can be computed with
{\schub}.
\end{proof}

By Remark (\ref{obs1}) we have a rational map
$\psi_{k}$  as in the diagram
\be\label{psik}  
\xymatrix
{
  \psi_{k}: \p{}(\D_{k}) \ar@{-->}[r]
\ar@{<-^)}[d]  &\p{}(\E_{k})\ar@{<-_)}[d]\\
              \psi_{k}^{-1}(\V_{k})
\vphantom{I^{I^I}}
        & 
\V_{k} \vphantom{I^{I^I}}
}
\ee
Set
$\Gamma_{k}:=\overline{\psi_{k}^{-1}(\V_{k})}\subset 
\p{}(\D_{k})$. 
Thus 
$$\Gamma_{k}=\{(p,[\w])\mid
\w_{p}=l^{k-1}(ydx-xdy)+h.o.t\, \text{ for
  some }\, l\in \Omega_{p}\}.
$$
The image
\\\centerline{$
C_{k}:=q(\Gamma_{k})\subset  \p N
$}  
parameterizes the space of foliations 
with the {\em MC} property.  

\bigskip
\coi{Lemma} \em
We have
\label{indet}
\\{\rm(i)} $\cod_{\p{}(\D_{k})}\Gamma_{k}
=\cod_{\p{}(\E_{k})} \V_{k}=k-2$. 
\\{\rm(ii)}
 Let $Z\subset \p{}(\D_{k})$ denote the
 indeterminacy locus of 
  $\psi_{k}$.
 Then  \begin{equation*}
\cod_{\p{}(\D_{k})}(Z)=k.
\end{equation*}  
\eco

\pf
Since
 $\psi_{k}$ (cf.\,diagram\,\ref{psik})
is induced by a surjective map of vector
bundles, its fibers
have the same dimension $n$. Therefore 
$\dim \Gamma_k=\dim \V_{k}+n$,
  and $\dim \p{}(\D_{k})=\dim
  \p{}(\E_{k})+n$.  Hence the
  equality for the codimension follows:
$\cod_{\p{}(\D_{k})}(Z)=\rk(\s{k-1}\Omega
\otimes\wedge^2\Omega (d+2))=k$.
\qed

\medskip
We may now find the degree of of the locus
of dicritical foliations with maximal contact.

\begin{coisa}\label{degC_{k}}
{\bf Proposition.} {\rm(i)}
The codimension of  $C_{k}$ in $\p N$ is 
$$\cod_{\p N}C_{k}=k^2+3k-2.
$$
{\rm(ii)}
The degree of $C_{k}$ is given by the  formula
$$\ba c
(k-1)\frac{1}{2}\left[\vphantom{\ba c1\\-\ea}
\frac{1}{4}(4k^6+20k^5-15k^4-66k^3+211k^2-218k+112)\right.
\\\left.
-(2k^5+7k^4+2k^3+24k^2-49k+44)d+(k^4+2k^3+10k^2+k+16)d^2\vphantom{\ba c1\\-\ea}\right].
\ea$$

\end{coisa}

\begin{proof} 

First of all, the restriction   
$q_{|\Gamma_{k}}$ is generically injective. 
For instance, it can be checked that 
the 1-form 
$$\w=(z_2^{d-(k-1)}(z_{0}+z_{1})^{k-1}+z_{0}^d+z_{1}^d)(z_{1}dz_{0}-z_{0}dz_{1})$$
has $p=[0:0:1]$ as its unique
singularity with order $k$ and $p$ is a reduced
point of the fiber 
$(q_{|\Gamma_{k}})^{-1}([\w])$.
To compute the codimension observe that
\\\centerline{$
\ba{rl}
\cod_{\p N}C_{k}=&\cod_{\p N}D_{k}+
\cod_{\p{}(\D_{k})}C_{k}=\\
& 
k(k+2)+k-2=k^2+3k-2.
\ea
$}
\\
Put $n=\dim \Gamma_{k}$.
By  lemma\,\ref{indet}\,(ii) we have that $\dim Z<n$,
hence 
$A_{n}(Z)=0$. Using the excision 
exact sequence 
(cf.\,\cite[Prop.\,1.8,\, p.\,21]{Fulton})
$$A_{n}(Z)\to A_{n}(\p{}(\D_{k}))\to A_{n}(\p{}(\D_{k})\setminus
Z)\to 0,$$ 
we deduce that 
$$
A_{n}(\p{}(\D_{k}))\simeq
A_{n}(\p{}(\D_{k})\setminus Z).
$$
Therefore, using that the class
$\psi_{k}^*[\V_{k}]$ is known  in
$\p{}(\D_{k})\setminus 
Z$, we can do the computations in 
$A_{n}(\p{}(\D_{k}))$.
Recalling (\ref{psik}) $\psi_{k}$
is a linear projection, we have
that 
$\psi_{k}^*H_{k}=H:=c_{1}(\OX_{\p{}(\D_{k})}(1))$.
Therefore
\begin{align*}
&\deg C_{k}=\deg q_{*}\Gamma_{k}=
\int H^n\cap [\Gamma_{k}]=
\int H^n\cap \psi_{k}^*[\V_{k}]= \\
&\int (u H^{n+k-2}+v hH^{n+k-3}+w h^2H^{n+k-4})\cap
     [\p{}(\D_{k})]= \tag{Lemma \ref{alpha}} \\
&\int (u H^{r+1}+v hH^{r}+w
h^2H^{r-1}) \cap [\p{}(\D_{k})],
\end{align*}
where $r=\rk \D_{k}$. Applying $p_{1*}$  and
using the  definition of Segre class we see
that what we are calculating is:
$$\int (u s_{2}(\D_{k})+v
s_{1}(\D_{k})+w s_{0}(\D_{k})\cap
[\p 2]. 
$$
From Lemma \ref{alpha} we know the values of $u,v,
w$. The classes  
$s_{1}(\D_{k})$ and $s_{2}(\D_{k})$ are known
from 
(\ref{defDk})  and the beginning of the 
proof of \ref{degMk}. We finish
using {\schub}.
\end{proof}

\section{Concluding remarks}
It is worth mentioning
 that for foliations of degree $d\geq2$,
the scheme of singularities completely
determines the foliation. Moreover, the schemes
of $d^2+d+1$ points that can occur as singular
scheme of a foliation are known,
cf.\,\cite{campillo}.
It would be nice  to work out the enumerative
geometry of the loci of foliations with
scheme of singularities subject to collisions.

The reader is invited to check that formulas
similar to  \ref{cor2.1}, \ref{degDk} and
\ref{degC_{k}} can be written down for an
arbitrary surface. Precisely, given a smooth,
projective surface $X$, we may fix an
ample divisor class $h$ and look at the space
of foliations $\p
N=\p{}(H^0(\Omega^1_X\otimes\OX((d+2)h)))$ for
$d\gg0$. The degree of $M_k$ can be written as
\\\centerline{$\ba c
\frac1{72}k(k+1)
(4k^4+8k^3-k^2-5k-6)c_1^2+
\left((24k^3+36k^2-12k-12)d+48k^3\right.
\\\left.
+72k^2-24k-24\right)c_1h+
6(k^2+6+4)c_2+
36(d+2)^2(k^2+k-1)h^2,
\ea$}
where we set for short $c_i=c_i\Omega^1_X$.
Substituting in the Chern numbers for \p2,
$(h^2=1,c_1h=-3,c_2=3)$\, reproduces
\ref{cor2.1}.
We include a script below.

\section{\sc schubert/maple script}
\label{script}
{\small
\begin{verbatim}
with(schubert): DIM:=2; omega:=bundle(2,c); f:=expand(Symm(k,omega)); 
g:=convert(%,list); s0:=sum(1, 'j'=0..k-1): s0:=factor(%);
s1:=sum('j', 'j'=0..k-1): s1:=factor(%); s2:=sum('j^2', 'j'=0..k-1): 
s2:=factor(%); s3:=sum('j^3', 'j'=0..k-1): s3:=factor(%);
G:=g;l:=[ ]: for i to nops(g) do  if has(g[i],k^3)then print(i):
l:=[op(l),i]: g[i]:=subs(k^3=s3,g[i]) fi od;g;
G:=g;l;for i to nops(g) do  if not i in l then
if has(g[i],k^2)then print(i):l:=[op(l),i]: g[i]:=subs(k^2=s2,g[i])
fi fi od;g;  G:=g;l;for i to nops(g) do
if not i in l then  if has(g[i],k) then print(i):
l:=[op(l),i]: g[i]:=subs(k=s1,g[i]) fi fi od; g;l; g[2]:=s0;
collect(convert(g,`+`),t); omega*o((d+2)*h); mtaylor(%%*%,t,3);
chern(2,%); factor(%);  #P2:c1^2=9*h^2,c2=3*h^2,c1=-3*h2
subs(c1^2=9*h^2,%); subs(c2=3*h^2,%); subs(c1=-3*h,%); print(indets(%));
factor(%); subs(h=1,%); collect(%,d);
\end{verbatim}
}

\end{document}